\documentclass[10pt,reqno]{amsart}
\usepackage{amsmath}
\usepackage{amssymb}
\usepackage{amsthm}

\newlength{\defbaselineskip}
\newcommand{\setlinespacing}[1]%
           {\setlength{\baselineskip}{#1 \defbaselineskip}}

\numberwithin{equation}{section}

\newtheorem{thm}{Theorem}[section]

\newtheorem{lem}[thm]{Lemma}

\theoremstyle{definition}

\theoremstyle{remark}
\newtheorem{rem}[thm]{Remark}
\numberwithin{equation}{section}

\begin{document}
\title[Dirac equations in modulation spaces]
{On Dirac equations with Hartree type nonlinearity in modulation spaces}

\author{Seongyeon Kim, Hyeongjin Lee and Ihyeok Seo}

\thanks{This research was supported by the POSCO Science Fellowship of POSCO TJ Park Foundation and  the Research Grant of Jeonju University in 2024 (S. Kim), and by NRF-2022R1A2C1011312 (I. Seo).}

\subjclass[2020]{Primary: 35Q40; Secondary: 42B35, 35A01 }
\keywords{Dirac equation, Modulation spaces, well-posedness}

\address{Department of Mathematics Education, Jeonju University, Jeonju 55069, Republic of Korea}
\email{sy\_kim@jj.ac.kr}

\address{Department of Mathematics, Sungkyunkwan University, Suwon 16419, Republic of Korea}
\email{hjinlee@skku.edu}
\email{ihseo@skku.edu}

\begin{abstract}
We obtain the local well-posedness for Dirac equations with a Hartree type nonlinearity derived by decoupling the Dirac-Klein-Gordon system. We extend the function space of initial data, enabling us to handle initial data that were not addressed in previous studies.
\end{abstract}

\maketitle

\section{Introduction}
We consider the Dirac equation with Hartree type nonlinearity
\begin{equation}\label{DE}
    \begin{cases}
    (-i\partial_t -i\alpha\cdot\nabla+m\beta) \psi = (\lambda|\cdot|^{-\gamma}\ast \langle \psi,\beta\psi\rangle)\beta\psi,\\
    \psi(0,x)=\psi_0 (x),
    \end{cases}
\end{equation}
where $\psi:\mathbb{R}^{1+d}\rightarrow \mathbb{C}^{n}$, $n=2^{\lfloor\frac{d+1}{2}\rfloor}$, is a spinor field regraded as a column vector, $m\ge0$ is mass, and $\alpha=(\alpha_1, \dots, \alpha_d),\beta$ are $n\times n$ Hermitian matrices satisfying the usual anticommutation relations:  
\begin{equation*}
\beta^{2}=I_{n}, \quad \alpha_{j}\beta+\beta\alpha_{j}=0, \quad \alpha_{j}\alpha_{k}+\alpha_{k}\alpha_{j}=2\delta_{jk}I_{n}
\end{equation*}
for all $1\le j,k \le d$. For $d=3$, the standard choice for such matrices is the Dirac matrices
\begin{equation*}
\alpha_{j}=
\begin{pmatrix}
0 & \sigma_{j}\\
\sigma_{j} & 0
\end{pmatrix},
\quad
\beta=
\begin{pmatrix}
I_{2} & 0\\
0 & -I_{2}
\end{pmatrix},
\end{equation*}
with the Pauli matrices
\begin{equation*}
\sigma_{1}=
\begin{pmatrix}
0 & 1\\
1 & 0
\end{pmatrix},
\quad
\sigma_{2}=
\begin{pmatrix}
0 & -i\\
i & 0
\end{pmatrix},
\quad
\sigma_{3}=
\begin{pmatrix}
1 & 0\\
0 & -1
\end{pmatrix}.
\end{equation*}
The potential $\lambda|x|^{-\gamma}$ here is of Coulomb type where $0<\gamma<d$ and $\lambda\not =0$ is a real constant.

When $d=2$ and $\gamma=1$, the equation \eqref{DE} can be viewed as a simplified model of Chern-Simons-Dirac system in Coulomb gauge \cite{BCM}.
When $d\ge3$ and $\gamma=d-2$, 
it is derived by decoupling the Dirac-Klein-Gordon system 
\begin{equation}\label{DKG}
    \begin{cases}
     (-i\partial_t-i\alpha\cdot\nabla + m \beta) \psi = \phi\beta\psi\\
     (\partial_t^2 - \Delta + M^2) \phi = \langle \psi , \beta \psi \rangle
    \end{cases}
\end{equation}
which arises as a model for the description of particle interactions in relativistic quantum mechanics (see \cite{BD} for more details).
Indeed, assume that a scalar field $\phi(t,x)=e^{\pm iMt}f(x)$ is a standing wave. Then the Klein-Gordon part of \eqref{DKG} becomes
\begin{equation*}
    -\Delta \phi = \langle \psi,\beta\psi\rangle
\end{equation*}
whose solution is 
\begin{equation*}
    \phi=\lambda|\cdot|^{-d+2}\ast \langle \psi,\beta\psi\rangle
\end{equation*} 
Putting this into the Dirac part of \eqref{DKG} yields the desired equation. 

Equation \eqref{DE} with a quadratic term $|\psi|^2$ replacing $\langle \psi,\beta\psi\rangle$ was derived by Chadam and Glassey \cite{CG} by uncoupling the Maxwell-Dirac system with a vanishing magnetic field.  

For equation \eqref{DE}, Cho, Lee and Ozawa \cite{CLO} proved the local well-posedness in $H^{s}(\mathbb{R}^2)$ for $s>\gamma-1$ with $1<\gamma<2$ in massive case $m>0$.
This was improved by Lee \cite{L} to $s>(\gamma-1)/2$ with $1\le \gamma < d$ ($d=2,3$) in both massive and massless cases $m\ge0$. 
He also showed the local well-posedness in $L^2(\mathbb{R}^2)$ for $1/2<\gamma<1$,
and in $H^{3/8}(\mathbb{R}^2)$ and $H^{1/4}(\mathbb{R}^3)$ for the remaining cases where $0<\gamma\le1/2$ and $0<\gamma<1$, respectively. 
In case of \eqref{DE} with a quadratic term $|\psi|^2$ replacing $\langle \psi,\beta\psi\rangle$, 
Machihara and Tsutaya \cite{MT} proved the local well-posedness in $H^{s}(\mathbb{R}^d)\,(d\ge 3)$ for $s>\frac{\gamma}{2d}+\frac 12$ with $2<\gamma<d$ in massive case. 
Since then, Nakamura and Tsutaya \cite{NT} extended the range of $\gamma$ to $3/2 < \gamma <d$ but for larger $s>\frac{\gamma}{2}+\frac1d$ $(d\ge 2)$. 
(See also \cite{H} for some related results.)

In this paper we obtain a local well-posedness in a wider function space containing $H^s(\mathbb{R}^d)$ for $s>\gamma/2$ with $0<\gamma< d$.
See Remark \ref{rem} below for details.
Compared with the results in \cite{L} mentioned above,
this range improves, when $d=2$, the range $s\ge 3/8$ ($0<\gamma\le1/2$) to $s>\gamma/2$. 
For $d=3$, it improves the range 
$s\ge1/4$ ($0<\gamma<1$) to $s>\gamma/2$ for
$0<\gamma<1/2$.
Furthermore, our range improves the aforementioned results in \cite{NT} for the quadratic case 
to the range $s>\gamma/2$ for all $0<\gamma<d$.

For the purpose, we utilize the modulation spaces whose definition is based on the short-time Fourier transform
\begin{equation*}
    V_gf(x,\xi)= \int_{\mathbb{R}^d} e^{- i y\cdot \xi}f(y)\overline{g(y-x)}dy
\end{equation*}
given a non-zero Schwartz function $g$;
for $1\le p,q\le \infty$ and $s\in\mathbb{R}$, the modulation space $M_s^{p,q}(\mathbb{R}^d,\mathbb{C}^n)$ is defined to be the Banach space of all tempered distributions 
$f\in\mathcal{S}'(\mathbb{R}^d,\mathbb{C}^n))$ such that
\begin{equation*}
    \|f\|_{M_{s}^{p,q}(\mathbb{R}^d,\mathbb{C}^n)}
=\bigg(\int_{\mathbb{R}^d}\Big(\int_{\mathbb{R}^d}|V_gf(x,\xi)|^p dx\Big)^{\frac{q}{p}} (1+|\xi|^2)^{\frac{sq}2} d\xi\bigg)^{\frac1q} < \infty,
\end{equation*}
for some non-zero smooth rapidly decreasing function $g\in\mathcal{S}(\mathbb{R}^d)$, with suitable modification for $p=\infty$ or $q=\infty$.
It seems that this norm depends on $g$, so it is worth noting that different choices of $g$ give equivalent norms. 
For simplicity, we denote $M^{p,q}=M_0^{p,q}$ when $s=0$.

The following is our main result.
\begin{thm}\label{main}
Let $d\ge1$ and $0<\gamma<d$. 
Then, 
there exists $T>0$ such that \eqref{DE}  has a unique solution $\psi\in C([0,T],X)$
for
\begin{equation*} 
X=
    \begin{cases}
    M^{p,q}(\mathbb{R}^d, \mathbb{C}^n) \quad \text{if} \quad 1\le p \le 2,\, 1\le q \le \frac{2d}{d+\gamma},\\
    M_{s}^{p,1}(\mathbb{R}^d,\mathbb{C}^n) \quad \text{if} \quad 1<p<\frac{d}{d-\gamma},\, s\ge0
    \end{cases}
\end{equation*}
if $\psi_0\in X$.
Moreover, if a maximal time $T$, say $T^\ast$, is finite, then 
\begin{equation}\label{maximal}
    \limsup_{t\rightarrow T^*} \|\psi(t,\cdot)\|_X = \infty.
\end{equation} 
\end{thm}

\begin{rem}\label{rem}
The argument in this paper can be also applied to equation \eqref{DE} with a quadratic term $|\psi|^2$ replacing $\langle \psi,\beta\psi\rangle$ and the same result holds.
For $s>\gamma/ 2$ with $0<\gamma<d$, we have the embedding $H^s\subset M^{2,\frac{2d}{d+\gamma}}\subset L^2$ (see Lemma \ref{emb} below). 
 Consequently, $X$ contains $H^s(\mathbb{R}^d)$ for $s>\gamma/2$ with $0<\gamma<d$.
Furthermore, it becomes possible to contain  functions not in $L^2$ because
$M^{p,q} \subset L^p$ for $1\le p<2$ and $1\le q\le\min\{p, \frac{2d}{d+\gamma}\}$.
In this regard, our results can handle initial data that were not addressed in previous studies. For more details, please refer to the paragraph immediately preceding the one introducing the definition of modulation spaces.
\end{rem}

The paper is organized as follows.
In Section \ref{sec2}, we list some basic properties of the modulation spaces used in the sequel. Section \ref{sec3} is devoted to proving Theorem \ref{main} by making use of a trilinear estimate (Lemma \ref{tri}) for the nonlinearity. 
This estimate is shown in the last section, Section \ref{sec4}.

Throughout this paper, the bracket notations $\langle \cdot , \cdot \rangle$ and $\langle \cdot, \cdot \rangle_{L^2}$ stand for complex inner product and $L^2$ inner product, respectively.
For a vector valued function $f=(f_1,\cdots,f_n)$, we define 
$|f|=\big(\sum_{j=1}^n|f_j|^2\big)^{\frac12}$, $\|f\|_{L^p}=\big\|\big(\sum_{j=1}^n|f_j|^2\big)^{\frac12}\big\|_{L^p},$ $\int_{\mathbb{R}^d} f(x)dx=\big(\int_{\mathbb{R}^d} f_1(x)dx, \cdots, \int_{\mathbb{R}^d} f_n(x)dx\big)$, and $\widehat f=(\widehat{f_1}, \cdots, \widehat{f_n})$.
We will also use the Fourier-Lebesgue spaces  $\mathcal{F}L^p(\mathbb{R}^d,\mathbb{C}^n)$ consisting of distributions $f\in \mathcal{S}' (\mathbb{R}^d,\mathbb{C}^n)$ such that
\begin{equation*}
    \|f\|_{\mathcal{F}L^p(\mathbb{R}^d,\mathbb{C}^n)}^p:= \|\widehat f \|_{L^{p}(\mathbb{R}^d,\mathbb{C}^n)}^p
    =\int_{\mathbb{R}^d} |\widehat f (\xi) |^p d\xi< \infty.
\end{equation*} 
The letter $C$ stands for a positive constant which may be different at each occurrence.
We finally denote $A\lesssim B$ to mean $A\leq CB$ with unspecified constants $C>0$.

\section{Preliminaries}\label{sec2}

We now list some basic properties of the modulation spaces.
Since $f=(f_1,...,f_n)\in M_{s}^{p,q}(\mathbb{R}^d,\mathbb{C}^{n})$ if and only if $f_j \in M_{s}^{p,q}(\mathbb{R}^d,\mathbb{C})$ for all $j$,
most of the properties extend to the vector-valued contexts. 

\begin{lem}[\cite{G, T, CKS}]\label{property}
Let $1\le p, q, p_i, q_i\le\infty$ and $s_i\in\mathbb{R}$.
\begin{enumerate}
    \item If $p_1\le p_2$, $q_1\le q_2$ and $s_1\ge s_2$, then $M_{s_1}^{p_1,q_1}(\mathbb{R}^d, \mathbb{C}^n) \hookrightarrow M_{s_2}^{p_2,q_2}(\mathbb{R}^d, \mathbb{C}^n)$.
    \item 
    $M^{p,q_1}(\mathbb{R}^d, \mathbb{C}^n) \hookrightarrow L^p(\mathbb{R}^d, \mathbb{C}^n) \hookrightarrow M^{p,q_2}(\mathbb{R}^d, \mathbb{C}^n)$ holds for $q_1\le \min\{p,p'\}$ and $q_2 \ge \max\{p,p'\}$.
    \item $M^{\min\{p',2\},p}(\mathbb{R}^d, \mathbb{C}^n)\hookrightarrow \mathcal{F}L^{p}(\mathbb{R}^d, \mathbb{C}^n)\hookrightarrow M^{\max\{p',2\},p}(\mathbb{R}^d, \mathbb{C}^n)$.
\end{enumerate}
\end{lem}

\begin{lem}[Theorem 2.8 and Remark 2.9 of \cite{Tr}]\label{multi}
Let $1\le p_i, q_i\le\infty$ and $s\ge0$.
Let $E_i$ be complex Banach spaces. If the map
\begin{equation*}
    \bullet : E_1 \times E_2 \rightarrow E_3,
    \quad
    (x_1,x_2)\mapsto x_3 = x_1 \bullet x_2
\end{equation*}
is a continuous bilinear operator with operator norm $\|\bullet\| \le 1$,
then
    \begin{equation*}
        M_{s}^{p_1,q_1}(\mathbb{R}^d,E_1) \bullet M_{s}^{p_2,q_2}(\mathbb{R}^d,E_2) \hookrightarrow M_{s}^{p_3,q_3}(\mathbb{R}^d,E_3)
    \end{equation*}
    for $\frac{1}{p_1}+\frac{1}{p_2}=\frac{1}{p_3}$ and $\frac{1}{q_1}+\frac{1}{q_2}=1+\frac{1}{q_3}$.
 Examples of the \textit{multiplication} $\bullet$ are
    the multiplication with scalars: $\mathbb{C} \times E \rightarrow E, (\lambda, x) \mapsto \lambda x$,
    and the complex inner product: $\mathbb{C}^n \times \mathbb{C}^n \rightarrow \mathbb{C}$, $ \langle z, w \rangle \mapsto \sum_{j=1}^{n} z_j \overline w_j$.
\end{lem}

The following proposition shows embedding relationships of the Sobolev spaces $W^{s,p}(\mathbb{R}^d,\mathbb{C}^n)$, defined by
$$\|f\|_{W^{s,p}}=\big\|\big((1+|\cdot|^2)^{s/2}\widehat f(\cdot)\big)^{\vee}\big\|_{L^p},$$
and the modulation spaces.

\begin{lem}[Theorem 3.8 of \cite{RSW}, p. 272]\label{emb}
Let $1\le p,q \le \infty$ and $s_1, s_2 \in \mathbb{R}$.
Then 
$$W^{s_1,p}(\mathbb{R}^d, \mathbb{C}^n) \subset M_{s_2}^{p,q}(\mathbb{R}^d, \mathbb{C}^n)$$ 
if one of the following conditions is satisfied:
    \begin{itemize}
     \item $q\ge p > 1$,\, $s_1 \ge s_2 + \tau(p,q)$;
     \item $p>q$,\, $s_1> s_2 + \tau(p,q)$;
    \end{itemize}
    and $M_{s_1}^{p,q}(\mathbb{R}^d, \mathbb{C}^n) \subset W^{s_2,p}(\mathbb{R}^d, \mathbb{C}^n)$ if
    one of the following conditions is satisfied:
    \begin{itemize}
     \item $q\le p <\infty$,\, $s_1 \ge s_2 + \sigma(p,q)$;
     \item $p<q$,\, $s_1> s_2 + \sigma(p,q)$;
    \end{itemize}
where
\begin{align*}
    \tau(p,q)=\max\bigg\{0,\,d(\frac{1}{q}-\frac{1}{p}),\, d(\frac{1}{q}+\frac1p-1)   \bigg\},\\
    \sigma(p,q)=\max\bigg\{0,\,d(\frac{1}{p}-\frac{1}{q}),\, d(1-\frac{1}{p}-\frac1q)   \bigg\}.
\end{align*}
\end{lem}

Since the Hartree type nonlinearity
involves convolutions with the Coulomb type potential $|x|^{-\gamma}$, the use of the Hardy-Littlewood-Sobolev inequality for fractional integrals $I_\gamma f=|\cdot|^{-\gamma}\ast f$ will be pivotal in proving Theorem \ref{main} in the following section.

\begin{lem}[Hardy-Littlewood-Sobolev inequality, \cite{St}]\label{hls}
Let $0<\gamma<d$, $1<p<q<\infty$ and
$1/q + 1 = 1/p + \gamma/ d$.
Then 
\begin{equation*}
    \| I_{\gamma} f \|_{L^q(\mathbb{R}^d, \mathbb{C})} \le C_{d,\gamma,p} \| f \|_{L^p (\mathbb{R}^d, \mathbb{C})}.
\end{equation*}
\end{lem}

The following is an analog of the above proposition for modulation spaces.

\begin{lem}[Proposition 3.2 of \cite{B3}]\label{hlsmod}
Let $0<\gamma<d$, $1<p_1<p_2<\infty$,
$1\le q \le \infty$, and $s\ge0$.
Then,
\begin{equation*}
    \| I_{\gamma} f\|_{M_s^{p_2,q}(\mathbb{R}^d, \mathbb{C})} \lesssim \|f\|_{M_s^{p_1,q}(\mathbb{R}^d, \mathbb{C})}    
\end{equation*}
if
\begin{equation*}
    \frac1 {p_2} +1 = \frac 1{p_1} +\frac\gamma d.
\end{equation*}
\end{lem}

\section{Proof of Theorem \ref{main}}\label{sec3}
In this section we prove Theorem \ref{main}, temporarily assuming some trilinear estimates (Lemma \ref{tri}) on the modulation spaces for the nonlinearity.

By Duhamel's principle, the solution to \eqref{DE} can be written in the form
\begin{equation}\label{sol}
U(t)\psi_0 - i \mathcal{N}\psi(t)
\end{equation}
where $U(t)$ is an operator-valued Fourier multiplier given by
\begin{equation*}
U(t,\xi)=e^{- i t (m\beta+\sum_{j=1}^{d} \alpha_j \xi_j)},
\end{equation*}
and
\begin{equation*}
    \mathcal{N}\psi(t)
    =\int_0^t U(t-s)
    \big[(\lambda|\cdot|^{-\gamma}\ast \langle \psi,\beta\psi\rangle )\beta\psi\big](s)ds.
\end{equation*}
By the contraction mapping principle, it suffices to 
show that the mapping $\Phi$ defined by $\Phi(\psi)(t)=U(t)\psi_0 - i \mathcal{N}\psi(t)$ is a contraction on 
\begin{equation*}
    B_{M,T} = \{ \psi \in C([0,T], X) : \|\psi\|_{C([0,T], X)}\le M \}
\end{equation*}
for appropriate values of $T,M>0$. 
To do so, we use the fixed time estimate for the multiplier $U(t)$, established in \cite{Tr}, and obtain a trilinear estimate for the nonlinearity in the next section.

\begin{lem}[Theorem 1.1 of \cite{Tr}]\label{fix}
Let $1\le p,q \le \infty$ and $s\in\mathbb{R}$. Then, 
\begin{equation}\label{fixest}
    \| U(t)\psi_0\|_{M_{s}^{p,q}} \lesssim (1+|t|)^{d|1/2-1/p|} \| \psi_0 \|_{M_{s}^{p,q}}
\end{equation}
for any $t\in\mathbb{R}$.
\end{lem}

\begin{lem}\label{tri}
Let $0<\gamma<d$ and 
\begin{equation*}
X=
\begin{cases}
M^{p,q}(\mathbb{R}^d) \quad \text{if} \quad 1\le p \le 2,\, 1\le q \le \frac{2d}{d+\gamma},\\
M_{s}^{p,1}(\mathbb{R}^d) \quad \text{if} \quad 1<p<\frac{d}{d-\gamma},\, s\ge0.
\end{cases}
\end{equation*}
Then we have
\begin{equation}\label{triest}
    \|(|\cdot|^{-\gamma}\ast \langle \psi_1,\beta\psi_2\rangle)\beta\psi_3\|_{X}
    \lesssim \prod_{j=1}^{3}\|\psi_j\|_{X}.
\end{equation}
\end{lem}
Assuming for the moment Lemma \ref{tri}, we first show that $\Phi(\psi)\in B_{M,T}$ for $\psi\in B_{M,T}$.
We apply \eqref{fixest} to the homogeneous term in \eqref{sol} to see
\begin{equation*}
    \|U(t)\psi_0\|_{X} \le  C(1+|t|)^{d|1/2-1/p|} \|\psi_0\|_{X}
\end{equation*}
and apply \eqref{triest} to the Duhamel term in \eqref{sol} to see
\begin{align*}
    \|\mathcal{N}\psi\|_{X} &\le
     C(1+|t|)^{d|1/2-1/p|}\int_0^t \big\|\big[(|\cdot|^{-\gamma}\ast \langle \psi,\beta\psi\rangle)\beta\psi\big](s)\big\|_{X}ds\\
    &\le
     C(1+|t|)^{d|1/2-1/p|}\int_0^{t} \|\psi(t)\|_{X}^3ds\\
    &\le
     CT(1+|t|)^{d|1/2-1/p|}\|\psi\|_{C([0,T],X)}^3,
\end{align*}
under the conditions in Lemma \ref{tri}. 
Hence, 
\begin{align*}
    \|\Phi(\psi)\|_{C([0,T],X)}
    \le
    C_T \|\psi_0\|_{X} +  C_T T\|\psi\|_{C([0,T],{X})}^3
\end{align*}
where $C_T= C(1+T)^{d|1/2-1/p|}$.
Now, if we set $M=2C_T\|\psi_0\|_{X}$ and 
take $T,M>0$ small so that $C_TTM^2 \le 1/2$,
we conclude
\begin{align*}
    \|\Phi(\psi)\|_{C([0,T],X)}
    \le \frac M 2 + \frac M 2 = M.
\end{align*}
Thus, $\Phi(\psi)$ belongs to $B_{M,T}$.

Next, we show that $\Phi$ is a contraction on $B_{M,T}$. Note that
\begin{align}\label{cont}
    (|\cdot|^{-\gamma}\ast&\langle \psi,\beta\psi\rangle)\beta\psi
    -
    (|\cdot|^{-\gamma}\ast\langle \phi,\beta\phi\rangle)\beta\phi\nonumber\\
    =&\;
    (|\cdot|^{-\gamma}\ast\langle \psi,\beta\psi\rangle)\beta(\psi-\phi)
    +
    \big(|\cdot|^{-\gamma}\ast(\langle \psi,\beta\psi\rangle
    -\langle \phi, \beta\phi\rangle)\big)\beta\phi\nonumber\\
    =&\;(|\cdot|^{-\gamma}\ast\langle \psi,\beta\psi\rangle)\beta(\psi-\phi)\nonumber\\
    &\qquad+ (|\cdot|^{-\gamma}\ast\langle \psi,\beta(\psi-\phi)\rangle)\beta\phi\nonumber
    +(|\cdot|^{-\gamma}\ast\langle \psi-\phi,\beta\phi\rangle)\beta\phi.
\end{align}
Applying Lemma \ref{tri} yields then  
\begin{align*}
    \|(|\cdot|^{-\gamma}\ast\langle\psi,\beta\psi\rangle)\beta\psi
    &-
    (|\cdot|^{-\gamma}\ast\langle \phi,\beta\phi\rangle)\beta\phi
    \|_{X}\\
     \lesssim&\;
(\|\psi\|_{X}^{2}+\|\psi\|_{X}\|\phi\|_{X}+\|\phi\|_{X}^2)\|\psi-\phi\|_{X}.
\end{align*}
Taking this into account, and considering the previous arguments, we get
\begin{align*}
    \|\Phi(\psi)-\Phi(\phi)\|_{X} 
    \le&\,
    \Big\|\int_0^t U(t-s)(\mathcal{A}\psi - \mathcal{A}\phi)ds\Big\|_{X} \nonumber\\
    \le&\,
    C(1+|t|)^{d|1/2-1/p|}\int_0^t \|\mathcal{A}\psi -\mathcal{A}\phi\|_{X}ds \nonumber\\
    \leq&\,
    C_T \int_0^tC(\psi,\phi)\|\psi-\phi\|_{X}ds\nonumber\\
    \leq&\,
    TC_T M^2 \|\psi-\phi\|_{C([0,T],X)}, 
\end{align*}
where $\mathcal{A}\psi = (|\cdot|^{-\gamma}\ast \langle \psi , \beta\psi \rangle)\beta\psi$.
Thus, 
\begin{equation*}
    \|\Phi(\psi)-\Phi(\phi)\|_{C([0,T],X)} 
    \le
    \frac{1}{2}
    \|\psi-\phi\|_{C([0,T],X)}
\end{equation*}
for sufficiently small $T>0$.

Finally, to show \eqref{maximal}, 
first suppose $\limsup_{t\rightarrow T^*}\|\psi(t)\|_{X} < \infty$. Then, there exists a sequence $\{t_n\}_{n=1}^{\infty}$  such that 
\begin{equation*}
 t_n \rightarrow T^* \quad \text{as}\quad n\rightarrow \infty
 \quad \text{and} \quad \|\psi(t_n)\|_{X} \le M
\end{equation*}
for some $M>0$.
For each $n\in\mathbb{N}$, considering initial data as $\psi(t_n,x)$ instead of $\psi(0,x)$, we can find a local solution
\begin{equation*}
    \psi_n \in C([t_n,t_n+T(M)],X)
\end{equation*}
where $T(M)$ depending only on $M$. Thus, there is a large $N\in\mathbb{N}$ such that $T^*<t_n+T(M)$ for each $n>N$. This contradicts to maximality of $T^*$.

\section{Proof of Lemma \ref{tri}}\label{sec4}
In this section, we prove the trilinear estimate for the nonlinearity.
We first show the case $X=M^{p,q}(\mathbb{R}^d)$.
By Lemma \ref{multi}, we have
\begin{align} \label{a1}
\nonumber
\|(|\cdot|^{-\gamma}\ast \langle \psi_1,\beta\psi_2\rangle)\beta\psi_3\|_{M^{p,q}}
&\lesssim
    \||\cdot|^{-\gamma}\ast \langle \psi_1,\beta\psi_2\rangle\|_{M^{\infty,1}}
    \|
    \psi_{3}\|_{M^{p,q}}\\
    \nonumber
&\lesssim \|\psi_{3}\|_{M^{p,q}}       \sum_{k=1}^{n}
    \||\cdot|^{-\gamma}\ast (\psi_{1,k}\overline{\psi_{2,k}})\|_{M^{\infty,1}}\nonumber\\
    &\lesssim
    \|\psi_{3}\|_{M^{p,q}}\sum_{k=1}^{n}
    \||\cdot|^{-\gamma}\ast (\psi_{1,k}\overline{\psi_{2,k}})\|_{\mathcal{F}L^{1}},
\end{align}
where $\psi_{j,k}$ is the $k$-th component of vector-valued function $\psi_j$.
Here, for the last inequality, we used Lemma \ref{property}, (3).
Using H\"older's inequality and Lemma \ref{hls}, the summand in \eqref{a1} is bounded as 
\begin{align}\label{a}
\nonumber
\||\cdot|^{-\gamma}\ast (\psi_{1,k}\overline{\psi_{2,k}})\|_{\mathcal{F}L^{1}}
    &\leq
    \iint_{\mathbb{R}^{2d}}
    \frac{|{\psi_{1,k}^{\vee}}(\eta-\xi)||\widehat{\overline{\psi_{2,k}}}(\eta)|}{|\xi|^{d-\gamma}}d\eta d\xi\\
    \nonumber
    &=\big\langle I_{d-\gamma}|\psi_{1,k}^{\vee}|,|\widehat{\overline{\psi_{2,k}}}|\big\rangle_{L_\eta^2}\\
    \nonumber
    &\leq\|I_{d-\gamma}|\psi_{1,k}^{\vee}|\|_{L^{\frac{2d}{d-\gamma}}} \|\widehat{\overline{\psi_{2,k}}}\|_{L^{\frac{2d}{d+\gamma}}}\\
    &\lesssim\|\psi_{1,k}^{\vee}\|_{L^{\frac{2d}{d+\gamma}}}
\|\widehat{\overline{\psi_{2,k}}}\|_{L^{\frac{2d}{d+\gamma}}}
\end{align}
when
\begin{equation*}
0<\gamma<d.
\end{equation*}
Then, by combining \eqref{a1} and \eqref{a} and using (3) in Lemma \ref{property}, we have 
\begin{align*}
\nonumber
&\|(|\cdot|^{-\gamma}\ast \langle \psi_1,\beta\psi_2\rangle)\beta\psi_3\|_{M^{p,q}}\\
&\qquad\quad\lesssim
    \|\psi_{3}\|_{M^{p,q}}\sum_{k=1}^n    \|\psi_{1,k}\|_{M^{\min\{\frac{2d}{d-\gamma},2\},\frac{2d}{d+\gamma}}}
    \|\overline{\psi_{2,k}}\|_{M^{\min\{\frac{2d}{d-\gamma},2\},\frac{2d}{d+\gamma}}}
\nonumber\\
&\qquad\quad= \|\psi_{3}\|_{M^{p,q}}\sum_{k=1}^n\|\psi_{1,k}\|_{M^{2,\frac{2d}{d+\gamma}}}
    \|\overline{\psi_{2,k}}\|_{M^{2,\frac{2d}{d+\gamma}}},
\end{align*}
which implies the desired estimates 
\begin{equation*}
    \|
    (|\cdot|^{-\gamma}\ast \langle \psi_1,\beta\psi_2\rangle)\beta\psi_3\|_{M^{p,q}}
    \lesssim \prod_{j=1}^{3}\|\psi_j\|_{M^{p,q}}
\end{equation*}
for $1\le p \le 2$ and $1\le q \le \frac{2d}{d+\gamma}$, by using (1) in Lemma \ref{property}.

It remains to show the case $X=M_s^{p,1}$.
By Lemma \ref{multi} and (1) in Lemma \ref{property}, we see that, for $0<\gamma<d$,
\begin{align*}
     \| (|\cdot|^{-\gamma}\ast \langle \psi_1, \beta \psi_2 \rangle) \beta\psi_3 \|_{M_s^{p,1}}
     &\lesssim
     \|I_\gamma\langle \psi_1, \beta \psi_2 \rangle\|_{M_s^{\infty,1}}\|\beta\psi_3\|_{M_s^{p,1}}\\
     &\lesssim
     \|I_\gamma\langle \psi_1, \beta \psi_2\rangle\|_{M_s^{\frac{pd}{d+p(\gamma-d)},1}}\|\psi_3\|_{M_s^{p,1}}.
\end{align*}
Then, by Lemmas \ref{hlsmod} and \ref{multi}, we obtain
\begin{align*}
\|I_\gamma\langle \psi_1, \beta \psi_2 \rangle\|_{M_s^{\frac{pd}{d+p(\gamma-d)},1}} &\lesssim \|\langle \psi_1, \beta \psi_2 \rangle\|_{M_s^{p,1}}\\
&\lesssim \|\psi_1\|_{M_s^{2p,1}}\|\beta\psi_2\|_{M_s^{2p,1}}\lesssim
    \prod_{j=1}^2\|\psi_j\|_{M_s^{p,1}},
\end{align*}
which concludes the proof.

\end{document}